\documentclass[a4paper,11pt,leqno]{amsart}
\usepackage{charter}
\usepackage{newtxmath}
\usepackage{tikz}
\usepackage[T1]{fontenc}
\usepackage{amsthm,amsmath}
\usepackage{todonotes}
\usepackage[alphabetic,initials,nobysame]{amsrefs}
\allowdisplaybreaks
\usepackage{xcolor}
\definecolor{darkred}{rgb}{0.5,0,0} %
\definecolor{darkblue}{rgb}{0,0,0.5} %
\usepackage[colorlinks,linkcolor=darkred,citecolor=darkblue,urlcolor=darkblue,hypertexnames=true]{hyperref}
\usepackage{enumitem}
\usepackage{comment}

\newcommand{\rn}{{\mathbb R}^n}
\newcommand{\C}{{\mathbb C}}

\newcommand{\dist}{\operatorname{dist}}

\DeclareMathOperator{\Lip}{Lip}

\numberwithin{equation}{section}
\theoremstyle{plain}
\newtheorem{theorem}{\textbf{THEOREM}}[section]
\newtheorem{lemma}[theorem]{\textsc{Lemma}}

\newtheorem{corollary}[theorem]{\textsc{Corollary}}

\theoremstyle{definition}
\newtheorem{example}[theorem]{\textsc{Example}}
\newtheorem{Thm}{Theorem}

\newtheorem{definition}[theorem]{\textsc{Definition}}

\theoremstyle{remark}
\newtheorem{remark}[theorem]{Remark}

\newtheorem*{theorem*}{Theorem}

\def\leq{\leqslant}
\def\geq{\geqslant}

\title{A complete classification of isometries of quasihyperbolic domains}
\author{Abhishek Pandey}
\dedicatory{Dedicated to David A. Herron and Pekka Koskela for their pioneering contributions to quasihyperbolic geometry.}

\address{Department of Mathematics and Statistics,
	University of Jyv\"askyl\"a,
	P.O. Box 35, 40014 University of Jyv\"askyl\"a, Finland}
\email{pandey.a.pandey@jyu.fi}


\thanks{A. Pandey was supported by the Research Council of Finland,
	project number 360505.}

\begin{document}
	
	\begin{abstract}
		We classify all surjective quasihyperbolic isometries
		between arbitrary proper planar domains. Apart from Euclidean similarities, the orientation-preserving
		isometries are reciprocal maps on sectors or punctured
		planes and explicit complex-power maps on logarithmic
		spiral sectors. In particular, orientation-preserving
		planar quasihyperbolic isometries need not be M\"obius.
		As consequences, we obtain similarity rigidity when
		either domain has finite inradius, for multiply connected
		domains other than punctured planes, and for $C^1$
		domains other than half-planes. The $C^1$ conclusion
		strengthens H\"ast\"o's corresponding result for planar
		$C^3$ domains. Our planar classification, together with Richard's classification in dimensions $n\geq3$, completely solves the quasihyperbolic isometry problem between
		proper subdomains of $\mathbb R^n$ in every dimension $n\geq2$. As an application, we extend to dimension two
		the characterization of pointed Gromov--Hausdorff
		convergence of quasihyperbolic domains obtained
		by Herron, Richard, and Snipes.
	\end{abstract}

	\maketitle
	\tableofcontents
	
	\section{Introduction and main result}
	\subsection{Background and main result} For a domain $\Omega\subsetneq\rn$, $n\geq2$, let $$\delta_\Omega(x)=\dist(x,\partial\Omega):=\inf_{p\in \partial \Omega}|x-p| \quad \text{ and } \quad R_\Omega=\sup_{x\in\Omega}\delta_\Omega(x).$$
	The {\it quasihyperbolic metric} $k_{\Omega}$ in $\Omega$ is defined by the length element $|dx|/\delta_{\Omega}(x)$, where $|dx|$ denotes the Euclidean length element. Thus, $k_{\Omega}(x,y)=\inf_{\gamma}\ell_{k}(\gamma)$, where the infimum is taken over all rectifiable curves $\gamma$ in $\Omega$ joining $x$ and $y$, and $\ell_{k}(\gamma)=\int_{\gamma} |dx|/\delta_{\Omega}$
	is the {\it quasihyperbolic length} of a rectifiable curve $\gamma$ in $\Omega$. It is well-known that $(\Omega, k_{\Omega})$ is complete proper and geodesic metric space; see \cite{BHK01}*{Proposition 2.8}. The pair $(\Omega, k_{\Omega})$ is a \emph{quasihyperbolic domain.} 
	
	The quasihyperbolic metric was introduced by Gehring and Palka \cite{GP76} in their study of quasiconformally homogeneous domain. It played an important role in Jones's classical characterization of $\mathrm{BMO}$ extension
	domains \cite{Jon80}. Soon after the characterizations of uniform domains by Gehring and Osgood \cite{GO79}, the metric becomes a basic tool in quasiconformal analysis. It has found numerous applications in Sobolev and quasiconformal removability \cite{JS00}, conformal rigidity \cite{NY20}, and the uniformization and Gromov hyperbolicity of metric spaces \cites{BHK01,BB03,BH20}. For further background and applications, see Koskela's
	survey \cite{Kos98}. 
	
	\medskip
	Martin \cite{Mar85}*{Corollary~4.8 and Example~4.11} proved that quasihyperbolic geodesic segments in arbitrary proper Euclidean domains are $C^{1,1}$ when parametrized by Euclidean arclength, and showed that the regularity is sharp. Later, Martin and Osgood
	\cite{MO86}*{Corollary~3.6} used the notion of generalized
	Gaussian curvature in the sense of Heins \cite{Hei62}
	and proved that the quasihyperbolic metric in every
	proper planar domain has nonpositive generalized
	Gaussian curvature.
	
	In 2007, H\"ast\"o \cite{H07} observed that its intrinsic
	geometry had received comparatively little attention,
	partly because the density $\delta_\Omega^{-1}$ need not
	be differentiable everywhere. Subsequent advances include
	Herron's proof that the quasihyperbolic metric universal
	cover of a planar domain is a Hadamard space \cite{Her21},
	and the theorem of Ikonen and the author \cite{IP26} that quasihyperbolic domains are $\mathrm{CAT}(2)$, with the constant $2$ optimal in dimensions $n\geq3$. This resolved several conjectures of V\"ais\"al\"a \cite{Vai09} on the local structure of the quasihyperbolic metric and its geodesics.
	
	\subsubsection*{Quasihyperbolic isometry problem} A complementary question concerns the maps preserving this
	metric. For domains $\Omega,\Omega'\subsetneq\rn$, $n\geq2$,
	a surjection $f:\Omega\to\Omega'$ is a
	\emph{quasihyperbolic isometry} if $f:(\Omega,k_\Omega)\to(\Omega',k_{\Omega'})$ is an isometry, i.e.,
	\begin{equation}\label{eq:isometry}
		k_{\Omega'}(f(x),f(y))=k_\Omega(x,y)
		\quad(x,y\in\Omega).
	\end{equation}
	The \emph{quasihyperbolic isometry problem} asks for a complete characterization of these maps, allowing both the source and target domains to vary \cite{H07}. Related isometry problems have been studied for the Apollonian
	metric \cite{HI05}, the distance-ratio and Seittenranta metrics \cite{HIL06}, and the Kulkarni--Pinkall metric \cite{HIM08}; see also the survey \cite{HIMPS07}. 
	
	In this paper, we give a complete
	classification of surjective quasihyperbolic isometries
	between proper planar domains, describing both the maps
	and the domains on which they occur; see Theorem~\ref{thm:main}. Together with Richard's higher-dimensional classification
	(Theorem~\ref{Thm:higher-dimension-classification}),
	it completely solves the quasihyperbolic isometry problem
	for proper subdomains of $\mathbb R^n$ in every dimension $n\geq2$.
	
	We first introduce a family of domains that appears in the
	classification of planar quasihyperbolic isometries; see
	Theorem~\ref{thm:main}.
	
\begin{figure}[htbp]
		\centering
		\begin{tikzpicture}[
			x=0.15\linewidth,
			y=0.15\linewidth,
			font=\small,
			boundary/.style={black,line width=0.65pt},
			window/.style={black!45,densely dashed,line width=0.4pt}
			]
			\begin{scope}
				\fill[black!12]
				(0,0)--(1,0)
				arc[start angle=0,end angle=90,radius=1]
				--cycle;
				\draw[window] (0,0) circle[radius=1];
				\draw[boundary] (1,0)--(0,0)--(0,1);
				\fill (0,0) circle[radius=0.018];
				
				\node[below left,inner sep=2pt] at (0,0) {$0$};
				\node at (0,1.30) {(a) Ordinary sector};
				\node at (0,-1.28) {$D_{0,0,\pi/2}\cap B(0,1)$};
			\end{scope}
			
			\begin{scope}[xshift=0.33\linewidth]
				\fill[black!12]
				plot[domain=-7:0,samples=300,variable=\t]
				({exp(\t)*cos(\t r)},
				{exp(\t)*sin(\t r)})
				arc[start angle=0,end angle=90,radius=1]
				plot[domain=0:-7,samples=300,variable=\t]
				({exp(\t)*cos((\t+pi/2) r)},
				{exp(\t)*sin((\t+pi/2) r)})
				--(0,0)--cycle;
				
				\draw[window] (0,0) circle[radius=1];
				\foreach \beta in {0,90}{
					\draw[boundary]
					plot[domain=-7:0,samples=300,variable=\t]
					({exp(\t)*cos(deg(\t)+\beta)},
					{exp(\t)*sin(deg(\t)+\beta)});
				}
				\fill (0,0) circle[radius=0.018];
				
				\node[above left,inner sep=2pt] at (0,0) {$0$};
				\node at (0,1.30) {(b) Spiral sector};
				\node at (0,-1.28) {$D_{1,0,\pi/2}\cap B(0,1)$};
			\end{scope}
			
			\begin{scope}[xshift=0.66\linewidth]
				\fill[black!12] (0,0) circle[radius=1];
				\draw[window] (0,0) circle[radius=1];
				\draw[boundary]
				plot[domain=-7:0,samples=300,variable=\t]
				({exp(\t)*cos(\t r)},
				{exp(\t)*sin(\t r)});
				\fill (0,0) circle[radius=0.018];
				
				\node[above left,inner sep=2pt] at (0,0) {$0$};
				\node at (0,1.30) {(c) Spiral slit};
				\node at (0,-1.28) {$D_{1,0,2\pi}\cap B(0,1)$};
			\end{scope}
		\end{tikzpicture}
		\caption{Shaded portions of the domains inside the unit disk.}
		\label{fig:spiral-domains}
	\end{figure}
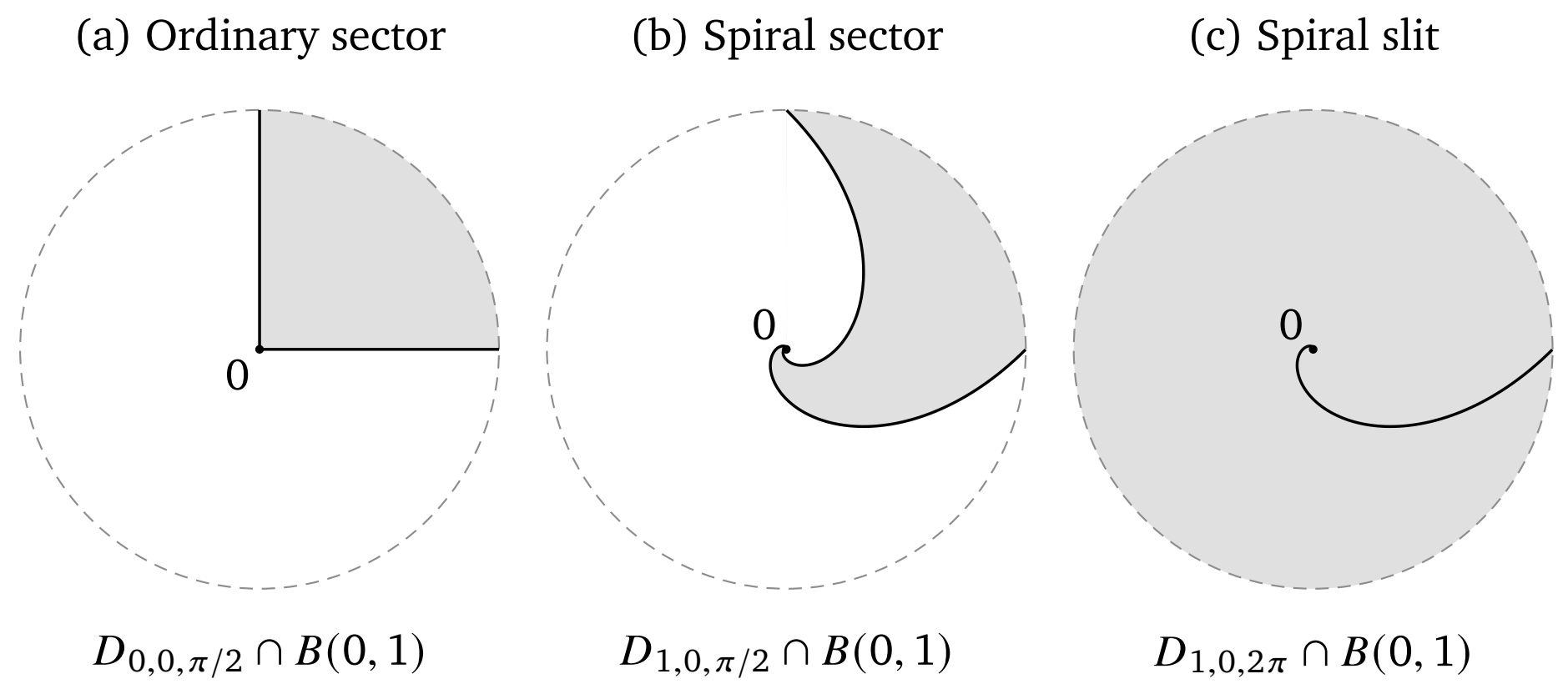

	\begin{example}\label{ex:main-example}
		For $\sigma,\theta\in\mathbb R$ and $0<\alpha\leq2\pi$, define
		\begin{equation}\label{eq:spiral-domains}
			\begin{split}
				U_{\sigma,\theta,\alpha}
				&=\{s+it\in\C:\theta<t-\sigma s<\theta+\alpha\},\\
				D_{\sigma,\theta,\alpha}
				&=\exp U_{\sigma,\theta,\alpha},\quad
				L_{\sigma,\theta,\alpha}
				=(\exp|_{U_{\sigma,\theta,\alpha}})^{-1}.
			\end{split}
		\end{equation}
		The map $\exp|_{U_{\sigma,\theta,\alpha}}$ is injective. Indeed if
		$w,w+2\pi im\in U_{\sigma,\theta,\alpha}$ for
		$m\in\mathbb Z$, then $2\pi|m|<\alpha\leq2\pi$, so $m=0$.
		Thus $L_{\sigma,\theta,\alpha}$ is a holomorphic logarithm on
		$D_{\sigma,\theta,\alpha}$.
		
		For $\sigma=0$, $D_{0,\theta,\alpha}$ is a sector,
		including a slit plane when $\alpha=2\pi$.
		For $\sigma\ne0$, $D_{\sigma,\theta,\alpha}$ is a logarithmic
		spiral sector. In both cases,
		\begin{equation}\label{eq:spiral-boundary}
			\partial D_{\sigma,\theta,\alpha}
			=\{0\}\cup
			\bigcup_{\beta\in\{\theta,\theta+\alpha\}}
			\{e^{(1+i\sigma)s+i\beta}:s\in\mathbb R\}.
		\end{equation}
		When $\alpha=2\pi$, the two boundary curves in
		\eqref{eq:spiral-boundary} coincide, since 
		$$ e^{(1+i\sigma)s+i(\theta+2\pi)}
		=e^{(1+i\sigma)s+i\theta},\quad s\in\mathbb R.$$ In this case,
		\begin{equation}\label{eq:spiral-slit}
			D_{\sigma,\theta,2\pi}
			=\C\setminus
			\left(\{0\}\cup
			\left\{e^{(1+i\sigma)s+i\theta}:s\in\mathbb R\right\}\right).
		\end{equation}
		For $\sigma\ne0$, \eqref{eq:spiral-slit} describes a plane slit
		along a logarithmic spiral; see Figure~\ref{fig:spiral-domains} for the illustration.
	\end{example}
	
	We now state our main result.
	\begin{theorem}\label{thm:main}
		Let $\Omega,\Omega'\subsetneq\C$ be domains and let
		$f:(\Omega,k_\Omega)\to(\Omega',k_{\Omega'})$ be a surjective
		isometry. After postcomposition with $z\mapsto\overline z$
		if necessary, exactly one of the following holds.
		\begin{enumerate}
			\renewcommand{\labelenumi}{(\arabic{enumi})}
			\renewcommand{\theenumi}{\thetheorem(\arabic{enumi})}
			
			\item\label{thm:main-similarity}
			\textbf{Similarity.}
			$f(z)=az+b$ for some $a\in\C\setminus\{0\}$ and $b\in\C$.
			
			\item\label{thm:main-reciprocal}
			\textbf{Reciprocal map.}
			For some $p\in\partial\Omega$, $A\in\C\setminus\{0\}$,
			and $B\in\C$,
			\begin{equation}\label{eq:main-reciprocal}
				f(z)=B+\frac{A}{z-p},\quad
				\Omega=p+D_{0,\theta,\alpha}\ \text{or}\ \C\setminus\{p\},
			\end{equation}
			where $\theta\in\mathbb R$ and $0<\alpha\leq2\pi$ in the sector case $\Omega=p+D_{0,\theta,\alpha}$.
			
			\item\label{thm:main-spiral}
			\textbf{Spiral-sector map.}
			For some $p\in\partial\Omega$, $A\in\C\setminus\{0\}$,
			$B\in\C$, $\sigma\in\mathbb R\setminus\{0\}$,
			$\theta\in\mathbb R$, and $0<\alpha\leq2\pi$,
			\begin{equation}\label{eq:main-power}
				\begin{split}
					\Omega&=p+D_{\sigma,\theta,\alpha},\quad
					\lambda_\sigma=\frac{-1+i\sigma}{1+i\sigma},\\
					f(z)&=B+A\exp\bigl(\lambda_\sigma
					L_{\sigma,\theta,\alpha}(z-p)\bigr).
				\end{split}
			\end{equation}
		\end{enumerate}
		Conversely, every map in above three cases is a
		quasihyperbolic isometry onto its image.
	\end{theorem}
	
	\begin{remark}\label{rem:non-Mobius-remark}
		The term
		$\exp\bigl(\lambda_\sigma
		L_{\sigma,\theta,\alpha}(z-p)\bigr)$
		in \eqref{eq:main-power} is the complex-power map since it is a complex power of $z-p$,
		defined using the holomorphic logarithm
		$L_{\sigma,\theta,\alpha}$. Writing
		$$
		(z-p)^{\lambda_\sigma}
		:=\exp\bigl(\lambda_\sigma
		L_{\sigma,\theta,\alpha}(z-p)\bigr),
		$$
		the maps in \eqref{eq:main-power} take the form
		$f(z)=B+A(z-p)^{\lambda_\sigma}$.
		
		For a locally univalent holomorphic map $f$, its Schwarzian
		derivative is
		\begin{equation}\label{eq:schwarzian-definition}
			S_f=\frac{f'''}{f'}-\frac32\left(\frac{f''}{f'}\right)^2
			=\left(\frac{f''}{f'}\right)'
			-\frac12\left(\frac{f''}{f'}\right)^2.
		\end{equation}
		Recall that a locally univalent holomorphic function on a
		domain is the restriction of a M\"obius transformation if
		and only if its Schwarzian derivative vanishes identically
		\cite{Osg98}*{pp.~275--276}. The maps in \eqref{eq:main-power} satisfy
		\begin{equation}\label{eq:main-schwarzian}
			|\lambda_\sigma|=1,\quad
			S_f(z)=\frac{1-\lambda_\sigma^2}{2(z-p)^2}\ne0.
		\end{equation}
		Hence they are not M\"obius. See Example~\ref{ex:spiral} for explicit computation with exponent $\lambda_1=i$. To the best of the author's knowledge, these are the first
		examples of orientation-preserving quasihyperbolic isometries between proper planar domains that are not restrictions of M\"obius transformations.
	\end{remark}
	
	It is known that every Euclidean similarity
	$$
	S(x)=cQx+b,\quad c>0,\quad b\in\rn,\quad Q\in O(n),
	$$
	satisfies \eqref{eq:isometry} with $\Omega'=S(\Omega)$.
	Here $O(n):=\{Q\in\mathbb{R}^{n\times n}:Q^TQ=I_n\}$ is the group of orthogonal linear automorphisms of $\rn$;
	equivalently, it is the group of Euclidean isometries of
	$\rn$ that fix the origin. Each $Q\in O(n)$ extends to
	$\rn\cup\{\infty\}$ by setting $Q(\infty)=\infty$. Examples of Euclidean similarities include translations,
	rotations, reflections, and dilations. Throughout the paper, M\"obius transformations of
	$\widehat{\mathbb R}^n=\mathbb R^n\cup\{\infty\}$
	are understood to include both orientations.
	In dimension two, these are the M\"obius and
	anti-M\"obius transformations in the usual
	complex-analytic terminology.
	
	The function $f:\mathbb H\to\mathbb H$ defined by
	$f(z)=-1/z$ is a M\"obius quasihyperbolic self-isometry
	of $\mathbb H=\{z\in\C:\operatorname{Im}z>0\}$,
	but is not a Euclidean similarity. More generally, in every dimension $n\geq2$, a M\"obius
	transformation of $\widehat{\mathbb R}^n$ restricts to
	a Euclidean similarity of $\mathbb R^n$ if and only
	if it fixes $\infty$. In particular, a planar M\"obius transformation
	$$
	M(z)=\frac{az+b}{cz+d}
	\quad\text{or}\quad
	M(z)=\frac{a\overline z+b}{c\overline z+d},
	\qquad a,b,c,d\in\C,\quad ad-bc\ne0,
	$$
	restricts to a Euclidean similarity of $\C$ if and only
	if $c=0$. Every M\"obius transformation that does not fix $\infty$ has the form
	\begin{equation}\label{eq:mobius-nonsimilarity-form}
		M(x)=b+cQ\frac{x-p}{|x-p|^2}
		\quad(x\in\rn\setminus\{p\}),
	\end{equation}
	where $p,b\in\rn$, $c>0$, and $Q\in O(n)$.
	Here $p=M^{-1}(\infty)$ and $b=M(\infty)$.
	
	\begin{definition}\label{def:similarity-rigidity}
		We say that \emph{similarity rigidity} holds for a domain
		$\Omega\subsetneq\mathbb R^n$, $n\ge 2$ if every quasihyperbolic isometry
		$f:(\Omega,k_\Omega)\to(f(\Omega),k_{f(\Omega)})$,
		where $f(\Omega)\subsetneq\mathbb R^n$ is a domain,
		is the restriction of a Euclidean similarity.
		Similarly, we say that \emph{M\"obius rigidity} holds
		for $\Omega$ if every such isometry is the restriction
		of a M\"obius transformation of $\widehat{\mathbb R}^n$.
		
		We say that $\Omega$ is \emph{self-similar rigid} if every
		surjective quasihyperbolic isometry
		$f:(\Omega,k_\Omega)\to(\Omega,k_\Omega)$
		is the restriction of a Euclidean similarity.
		Similarly, we say that $\Omega$ is
		\emph{self-M\"obius rigid} if every such self-isometry
		is the restriction of a M\"obius transformation of
		$\widehat{\mathbb R}^n$.
	\end{definition}
	
	\subsection{Known results on quasihyperbolic isometry}
	Let $\Omega,\Omega'\subsetneq\mathbb{R}^n$, $n\geq2$,
	be domains, and let
	$f:(\Omega,k_\Omega)\to(\Omega',k_{\Omega'})$
	be a surjective isometry. Martin and Osgood proved that,
	after postcomposition with a reflection if necessary,
	$f$ is conformal \cite{MO86}*{Theorem~2.6}.
	Consequently, if $n=2$, then $f$ is either holomorphic
	or antiholomorphic, and if $n\geq3$, then, by Liouville's theorem,
	$f$ is the restriction of a M\"obius transformation of
	$\widehat{\mathbb{R}}^n$.
	\subsubsection{H\"ast\"o's similarity rigidity results} We now review some of H\"ast\"o’s results from \cite{H07}. H\"ast\"o \cite{H07}*{p.~317} separates the isometry problem into conformality, M\"obius rigidity, and similarity rigidity.
	His Example~2.4 gives nonsimilarity isometries on sectors
	and punctured planes. Moreover, H\"ast\"o proved the following similarity rigidity results for
	\begin{itemize}
		\item planar $C^2$ domains other than half-planes,
		strictly convex domains, and strictly concave domains;
		see \cite{H07}*{Corollary~4.2}, where strict concavity
		means strict convexity of the complement;
		\item planar $C^3$ domains other than half-planes;
		see \cite{H07}*{Theorem~4.3}. His method uses
		generalized Gaussian curvature and its gradient;
	\end{itemize}
	
	For planar $C^1$ domains other than half-planes,
	H\"ast\"o also stated similarity rigidity under the
	additional assumption that the quasihyperbolic isometry
	is M\"obius \cite{H07}*{Proposition~2.2}, and deduced
	similarity rigidity for $C^1$ domains other than
	half-spaces in dimensions $n\geq3$
	\cite{H07}*{Corollary~2.3}. Richard pointed out a gap in this argument; see \cite{R19}*{Section~7.1, Example~7.1.2}, for a detailed discussion. Richard's inversion classification
	\cite{R19}*{Theorem~7.1.5}, which applies in every dimension
	$n\geq2$, yields an alternative proof of both conclusions.
	
	A domain $\Omega$ is a $C^m$ domain if its boundary is locally the graph of a $C^m$ function after a Euclidean change of coordinates.
	
	\subsubsection{Richard's higher-dimensional classification}
	
	Richard classifies quasihyperbolic isometries induced by
	spherical inversion in every dimension $n\geq2$
	\cite{R19}*{Theorem~7.1.5}. For $n\geq3$, conformality and
	Liouville's theorem reduce the general problem to that
	inversion case; see \cite{R19}*{Theorem~7.1.6}.
	We state a reformulation of these results which gives a
	complete classification in dimensions $n\ge 3$.
	For completeness, we include a short proof of this
	reformulation in our notation in Section~\ref{sec:higher}.
	
	\begin{Thm}\label{Thm:higher-dimension-classification}
		Let $n\geq3$, let $\Omega,\Omega'\subsetneq\rn$ be domains, and let $f:(\Omega,k_\Omega)\to(\Omega',k_{\Omega'})$ be a surjective map. Then $f$ is an isometry if and only if either $f$ is the
		restriction of a Euclidean similarity or there exist
		$p\in\partial\Omega$, $b\in\rn$, $c>0$, and $Q\in O(n)$
		such that
		\begin{equation}\label{eq:higher-classification}
			f(x)=b+cQ\frac{x-p}{|x-p|^2}\quad(x\in\Omega),\quad
			p+t(\Omega-p)=\Omega\quad(t>0).
		\end{equation}
		In the second case,
		\begin{equation}\label{eq:higher-image}
			\Omega'=b+cQ(\Omega-p).
		\end{equation}
	\end{Thm}
	
	\subsection{Consequences of the main result}
	The classification in Theorem~\ref{thm:main} gives the following similarity rigidity consequences, including similarity rigidity for domains of finite inradius and $C^1$ domains. Recall that
	$R_\Omega=\sup_{z\in\Omega}\delta_\Omega(z)$
	is the inradius of $\Omega$.
	
	\begin{corollary}\label{cor:consequences}
		Let $\Omega,\Omega'\subsetneq\C$ be domains, and let
		$f:\Omega\to\Omega'$ be a surjective quasihyperbolic isometry.
		Then the following statements hold.
		\begin{enumerate}[
			label=(\arabic*),
			ref=\thecorollary(\arabic*)
			]
			\item\label{cor:inradius}
			If $\min\{R_\Omega,R_{\Omega'}\}<\infty$, then $f$ is
			a Euclidean similarity. In particular, the conclusion
			holds if one of the domains is bounded.
			
			\item\label{cor:connectivity}
			If $\Omega$ is multiply connected and
			$\Omega\ne\C\setminus\{p\}$ for every $p\in\C$, then
			$f$ is a Euclidean similarity. Thus, a nonsimilarity
			quasihyperbolic isometry can occur only when $\Omega$
			is simply connected or a punctured plane.
			
			\item\label{cor:c1}
			If $\Omega$ is a $C^1$ domain other than a half-plane,
			then $f$ is a Euclidean similarity.
		\end{enumerate}
	\end{corollary}
	
	Corollary~\ref{cor:inradius} and
	Theorem~\ref{Thm:higher-dimension-classification}
	yield the following similarity rigidity result
	in every dimension $n\geq2$.
	
	\begin{corollary}\label{cor:inradius-all-dimensions}
		Let $n\geq2$, let $\Omega,\Omega'\subsetneq\rn$ be domains, and let
		$f:(\Omega,k_\Omega)\to(\Omega',k_{\Omega'})$ be a surjective
		isometry. If $\min\{R_\Omega,R_{\Omega'}\}<\infty$, then $f$ is the restriction of a Euclidean similarity.
	\end{corollary}
	
	The next consequence tells us that the quasihyperbolic metric determines the planar domain up to Euclidean similarity, even though an individual isometry can be non-Möbius.
	
	\begin{corollary}\label{cor:similar-domains}
		If $\Omega,\Omega'\subsetneq\C$ are quasihyperbolically isometric, then $\Omega'=S(\Omega)$ for a Euclidean similarity $S$, possibly reversing orientation.
	\end{corollary}
	
	Corollary~\ref{cor:similar-domains} allows us to test
	similarity rigidity and M\"obius rigidity using only
	quasihyperbolic self-isometries; see Definition~\ref{def:similarity-rigidity}
	
	\begin{corollary}\label{cor:self-rigidity}
		Let $\Omega\subsetneq\C$ be a domain. Then
		\begin{equation}\label{eq:self-similarity-rigidity}
			\Omega\text{ has similarity rigidity}
			\iff
			\Omega\text{ is self-similar rigid},
		\end{equation}
		and
		\begin{equation}\label{eq:self-mobius-rigidity}
			\Omega\text{ has M\"obius rigidity}
			\iff
			\Omega\text{ is self-M\"obius rigid}.
		\end{equation}
	\end{corollary}
	
	\begin{remark}
		By Theorem~\ref{thm:main} and
		Corollary~\ref{cor:self-rigidity},
		$\Omega$ is self-similar rigid if and only if
		$\Omega\notin\mathcal E$, where
		$$\mathcal E=\{\mathbb C\setminus\{p\}:p\in\mathbb C\}\cup
		\{p+D_{\sigma,\theta,\alpha}:
		p\in\mathbb C,\ \sigma,\theta\in\mathbb R,\
		0<\alpha\le2\pi\},$$
		and it is self-M\"obius rigid
		if and only if $\Omega\notin\mathcal N$, where
		$$\mathcal N=
		\{p+D_{\sigma,\theta,\alpha}:
		p\in\C,\ \sigma\in\mathbb R\setminus\{0\},\
		\theta\in\mathbb R,\ 0<\alpha\leq2\pi\}.$$
	\end{remark}
	
	We next explain how Theorem~\ref{thm:main}
	extends the known results.
	\begin{remark}\label{rem:non-similarities-qh-isometries}
		Theorem~\ref{thm:main} and its corollaries extend the known
		results as follows.
		\begin{itemize}
			\item Theorem~\ref{thm:main} gives a classification
			without boundary regularity assumptions.
			Corollary~\ref{cor:c1} strengthens H\"ast\"o's
			similarity rigidity for planar $C^3$ domains other
			than half-planes to the corresponding result for
			$C^1$ domains. It also removes the M\"obius assumption
			from \cite{H07}*{Proposition~2.2}.
			
			\item Theorem~\ref{thm:main} completely characterizes
			the orientation-preserving non-M\"obius
			quasihyperbolic isometries through the complex-power
			maps in \eqref{eq:main-power}.
			
			\item Contrary to the expectation expressed in
			\cite{H07}*{Example~2.4} and
			\cite{HIMPS07}*{Section~2}, sectors and punctured
			planes do not exhaust the domains admitting
			nonsimilarity quasihyperbolic isometries:
			the logarithmic spiral sectors in
			Theorem~\ref{thm:main} provide additional examples.
		\end{itemize}
	\end{remark}
	
	\subsection{Application to pointed Gromov--Hausdorff convergence}\label{subsec:pgh-application}
	Using the classification of Richard in dimensions $n\geq3$,
	Herron, Richard, and Snipes
	\cite{HRS20}*{Theorem~5.8 and Corollary~5.9}
	characterized pointed Gromov--Hausdorff convergence
	of quasihyperbolic domains in terms of uniform
	convergence of distance-to-complement functions on bounded
	Euclidean subsets of the limit domain. The normalized convergence result already appears in Richard's thesis
	\cite{R19}*{Theorem~7.2.1}, which also contains a
	corresponding planar result when the limit domain has
	$C^3$ boundary \cite{R19}*{Theorem~7.2.2}.
	The implication from boundary-distance convergence to
	pointed Gromov--Hausdorff convergence follows from the
	more general convergence theorem
	\cite{HRS20}*{Theorem~4.4}.
	Using Theorem~\ref{thm:main}, we establish the converse
	implication in dimension two without boundary regularity
	assumptions, thereby extending the characterization of
	Herron, Richard, and Snipes to arbitrary proper planar
	domains.
	
	For a proper domain $U\subsetneq\C$, define
	$d_U:\C\to[0,\infty)$ by
	\begin{equation}\label{eq:extended-boundary-distance}
		d_U(z)=\operatorname{dist}(z,\C\setminus U)
		\quad(z\in\C).
	\end{equation}
	Thus $d_U=\delta_U$ on $U$.
	We obtain the following characterization.
	
	\begin{theorem}\label{thm:planar-pgh}
		Let $\Omega,\Omega_j\subsetneq\C$ be domains, and let
		$a\in\Omega$ and $a_j\in\Omega_j$.
		The following are equivalent.
		\begin{enumerate}
			\renewcommand{\labelenumi}{(\arabic{enumi})}
			\renewcommand{\theenumi}{\thetheorem(\arabic{enumi})}
			
			\item\label{thm:planar-pgh-metric}
			The pointed quasihyperbolic spaces satisfy
			\begin{equation}\label{eq:planar-pgh-convergence}
				(\Omega_j,k_{\Omega_j};a_j)
				\xrightarrow{\mathrm{pGH}}
				(\Omega,k_\Omega;a).
			\end{equation}
			
			\item\label{thm:planar-pgh-distance}
			There exist $Q_j\in O(2)$ such that the similarities
			\begin{equation}\label{eq:planar-pgh-similarities}
				S_j(z)
				=a+\frac{\delta_\Omega(a)}{\delta_{\Omega_j}(a_j)}
				Q_j(z-a_j)
				\quad(z\in\C)
			\end{equation}
			satisfy
			\begin{equation}\label{eq:planar-pgh-uniform}
				\lim_{j\to\infty}
				\sup_{z\in\Omega\cap B(a,R)}
				\left|
				d_{S_j(\Omega_j)}(z)-\delta_\Omega(z)
				\right|=0
				\quad(R>0).
			\end{equation}
		\end{enumerate}
		The convergence in \eqref{eq:planar-pgh-uniform}
		is equivalent to
		$d_{S_j(\Omega_j)}(z)\to\delta_\Omega(z)$
		for every $z\in\Omega$.
	\end{theorem}
	
	\subsection{Idea of the proof and organisation of the paper}\label{subsec:idea-of-the-proof}
	We prove Theorem~\ref{thm:main} in Section~\ref{subsec:proof-of-main-result}. The idea of the proof is as follows.
	Let $f:\Omega\to\Omega'$ be a surjective quasihyperbolic
	isometry. After postcomposition with complex conjugation
	if necessary, we may assume that $f$ is holomorphic.
	Write $P=f''/f'$ and $u(z)=\log|f'(z)|$.
	We divide the proof into the cases $P\equiv0$ and
	$P\not\equiv0$. If $P\equiv0$, then $f(z)=az+b$
	with $a\ne0$.
	
	Suppose that $P\not\equiv0$. In Section~\ref{sec:differential-identity}, we use the
	boundary-distance identity
	\begin{equation}\label{eq:intro-density}
		\delta_{\Omega'}(f(z))=|f'(z)|\delta_\Omega(z)
		\quad(z\in\Omega)
	\end{equation}
	and the eikonal identities $|\nabla\delta_\Omega|=1$ and
	$|\nabla\delta_{\Omega'}|=1$ at their respective
	differentiability points to derive the key differential identity
	\begin{equation}\label{eq:intro-key-identity}
		\nabla\delta_\Omega\cdot\nabla u
		=-\frac{\delta_\Omega}{2}|\nabla u|^2.
	\end{equation}
	This holds at every differentiability point of
	$\delta_\Omega$. See Remark~\ref{rem:differential-identity-literature} for the relation with earlier work. Introducing $q=-2/P$ on $G=\{z\in\Omega:P(z)\ne0\}$ gives its vector form
	\begin{equation}\label{eq:intro-vector-identity}
		\nabla\delta_\Omega\cdot q=\delta_\Omega,
	\end{equation}
	which holds at every differentiability point of
	$\delta_\Omega$ in $G$.
	
	In Section~\ref{sec:classification}, we prove the key
	observation that restricting the vector identity \eqref{eq:intro-vector-identity} to nearest-boundary segments and differentiating yields
	$\operatorname{Re}q'=1$ on $G$.
	Holomorphicity and the identity theorem then imply
	that $P$ has no zeros in $\Omega$, so $G=\Omega$, and
	that there exist $\sigma\in\mathbb R$ and
	$p\in\C\setminus\Omega$ such that
	\begin{equation}\label{eq:intro-affine-vector}
		P(z)=-\frac{2}{(1+i\sigma)(z-p)},\quad
		q(z)=(1+i\sigma)(z-p)
		\quad(z\in\Omega).
	\end{equation}
	In particular, $q$ is affine. Since $G=\Omega$ and
	$\delta_\Omega$ is differentiable a.e. in $\Omega$,
	\eqref{eq:intro-vector-identity} and
	\eqref{eq:intro-affine-vector} give
	$$
	\nabla\delta_\Omega(z)\cdot(1+i\sigma)(z-p)
	=\delta_\Omega(z)
	\quad\text{for a.e. }z\in\Omega.
	$$
	By Lemma~\ref{lem:flow}, it follows that $p\in\partial\Omega$
	and that $\Omega$ is invariant under the flow
	$$
	\Phi_t(z)=p+e^{(1+i\sigma)t}(z-p),
	\quad t\in\mathbb R.
	$$
	When $\sigma=0$, this flow consists of dilations
	about $p$. When $\sigma\ne0$, it combines rotations
	and dilations about $p$, and its orbits in $\Omega$
	are logarithmic spirals. Connectedness then gives either
	$\Omega=\C\setminus\{p\}$ or
	$\Omega=p+D_{\sigma,\theta,\alpha}$ for some
	$\theta\in\mathbb R$ and $0<\alpha\leq2\pi$;
	see Lemma~\ref{lem:spiral-domains}.
	If $\sigma=0$, integrating the formula for $P=f''/f'$
	in \eqref{eq:intro-affine-vector} gives
	\eqref{eq:main-reciprocal}.
	If $\sigma\ne0$, a winding-number argument rules out
	$\Omega=\C\setminus\{p\}$, and integration gives
	\eqref{eq:main-power}.
	For the converse, we verify that the listed maps are
	conformal bijections onto their images and satisfy
	\eqref{eq:intro-density}.
	
	We prove Theorem~\ref{thm:planar-pgh} in
	Section~\ref{sec:pgh-proof}.
	
	\subsection*{Acknowledgements}
	The author thanks Kai Rajala for his comments on an earlier
	version of this manuscript, which helped improve the exposition and clarity of the manuscript. In particular, his question about whether the differential identity \eqref{eq:key-identity} had appeared previously prompted the literature search discussed in Remark~\ref{rem:differential-identity-literature}. The author thanks David A. Herron and Pekka Koskela for their encouragement and discussion.
	
	\section{The differential identity}\label{sec:differential-identity}
	
	For a planar domain $D$, write $ D^*=\{\overline z:z\in D\}$. Thus $D^*$ denotes reflection, whereas $\overline D$ denotes Euclidean closure. We identify vectors in $\mathbb R^2$ with complex numbers and use
	$v\cdot w=\operatorname{Re}(v\overline w)$.
	
	\begin{lemma}\label{lem:density}
		Let $\Omega,\Omega'\subsetneq\C$ be domains, and let
		$f:\Omega\to\Omega'$ be a holomorphic bijection.
		Then $f$ is a quasihyperbolic isometry if and only if
		\begin{equation}\label{eq:density}
			\delta_{\Omega'}(f(z))
			=|f'(z)|\delta_\Omega(z)
			\quad(z\in\Omega).
		\end{equation}
	\end{lemma}
	
	\begin{proof}
		Recall the infinitesimal formula in \cite{H07}*{p.~321}
		\begin{equation}\label{eq:infinitesimal}
			\lim_{w\to z}\frac{k_\Omega(z,w)}{|z-w|}
			=\frac1{\delta_\Omega(z)}.
		\end{equation}
		If $f$ is a quasihyperbolic isometry,
		\eqref{eq:infinitesimal} gives \eqref{eq:density}.
		Conversely, \eqref{eq:density} gives
		\begin{equation}\label{eq:length-preservation}
			\int_{f\circ\gamma}\frac{|dx|}{\delta_{\Omega'}}
			=\int_\gamma\frac{|f'||dx|}{\delta_{\Omega'}\circ f}
			=\int_\gamma\frac{|dx|}{\delta_\Omega}
		\end{equation}
		for every rectifiable curve $\gamma\subset\Omega$. Taking infimum in \eqref{eq:length-preservation} over all rectifiable curves joining $x$ and $y$, and
		applying the same argument to $f^{-1}$, gives
		\eqref{eq:isometry}.
	\end{proof}
	
	We recall the standard properties of $\delta_\Omega$;
	see \cite{MM03}*{Section~3}.
	The function $\delta_\Omega$ is $1$-Lipschitz and is
	differentiable at $z\in\Omega$ if and only if $z$ has a
	unique nearest boundary point. If $\xi\in\partial\Omega$
	is this point, then
	$\nabla\delta_\Omega(z)=(z-\xi)/|z-\xi|$. Consequently,
	\begin{equation}\label{eq:eikonal}
		\begin{gathered}
			\delta_\Omega\in\Lip(\Omega),\,\,\,
			|\nabla\delta_\Omega(z)|=1
			\quad\text{at every differentiability point }z\in\Omega.
		\end{gathered}
	\end{equation}
	By Rademacher's theorem, $\delta_\Omega$ is
	differentiable a.e. in $\Omega$. Our discussion is closely related to the geometry of the
	medial axis of $\Omega$. For related work on medial axis and the regularity of distance functions we refer to \cites{CF79,CCM97,Dam03}.
	
	The following lemma shows that $\delta_\Omega$ is
	differentiable at every interior point of a segment joining
	a point of $\Omega$ to one of its nearest boundary points,
	and that its gradient is constant along this segment.
	
	\begin{lemma}\label{lem:nearest-segment}
		Let $\Omega\subsetneq\rn$, $n\ge 2$
		be a domain. Let $x\in\Omega$, set $r=\delta_\Omega(x)$,
		and choose $\xi\in\partial\Omega$ with $|x-\xi|=r$.
		Put $\nu=(x-\xi)/r$.
		For every $0<t<r$, the point $z_t=\xi+t\nu$
		has $\xi$ as its unique nearest boundary point.
		Consequently,
		\begin{equation}\label{eq:nearest-segment}
			z_t\in\Omega,\quad
			\delta_\Omega(z_t)=t,\quad
			\nabla\delta_\Omega(z_t)=\nu
			\quad(0<t<r).
		\end{equation}
	\end{lemma}
	
	\begin{proof}
		Since $|x-z_t|=r-t<r$, we have $z_t\in B(x,r)\subset\Omega$.
		For every $\eta\in\partial\Omega$,
		\begin{equation}\label{eq:nearest-comparison}
			|z_t-\eta|\geq |x-\eta|-|x-z_t|
			\geq r-(r-t)=t.
		\end{equation}
		Equality in \eqref{eq:nearest-comparison} requires
		$|x-\eta|=r$ and equality in the triangle inequality.
		Thus $x-z_t$ and $z_t-\eta$ are nonnegative multiples of $\nu$,
		and $z_t-\eta=t\nu$. Hence $\eta=\xi$.
		The distance-function properties imply \eqref{eq:nearest-segment}.
	\end{proof}
	
	Now, let $f:\Omega\to\Omega'$
	be a holomorphic quasihyperbolic isometry and write
	\begin{equation}\label{eq:analytic-notation}
		u(z)=\log|f'(z)|,\quad P(z)=\frac{f''(z)}{f'(z)}.
	\end{equation}
	Since $f$ is injective and holomorphic,
	$f'(z)\ne0$ for every $z\in\Omega$.
	Thus $u$ and $P$ are globally defined, and
	\begin{equation}\label{eq:u}
		\Delta u=0,\quad P=u_x-iu_y\quad\text{in }\Omega.
	\end{equation}
	
	Comparing the eikonal identities \eqref{eq:eikonal} for
	$\delta_\Omega$ and $\delta_{\Omega'}$ through the
	boundary-distance identity \eqref{eq:density}, we obtain
	the following differential identity.
	
	\begin{lemma}\label{lem:identity}
		Let $\Omega,\Omega'\subsetneq\C$ be domains, and let
		$f:\Omega\to\Omega'$ be a surjective holomorphic
		quasihyperbolic isometry. Then at every differentiability point $z\in\Omega$
		of $\delta_\Omega$,
		\begin{equation}\label{eq:key-identity}
			\nabla\delta_\Omega(z)\cdot\nabla u(z)
			=-\frac{\delta_\Omega(z)}{2}|\nabla u(z)|^2.
		\end{equation}
		Let $G=\{z\in\Omega:P(z)\ne0\}$ and define
		$q=-2/P$ on $G$.
		Then, at every differentiability point $z\in G$
		of $\delta_\Omega$,
		\begin{equation}\label{eq:vector-identity}
			\nabla\delta_\Omega(z)\cdot q(z)=\delta_\Omega(z).
		\end{equation}
	\end{lemma}
	
	\begin{proof}
		Fix a differentiability point $z\in\Omega$ of $\delta_\Omega$.
		By \eqref{eq:density}, the identity
		$\delta_{\Omega'}=(e^u\delta_\Omega)\circ f^{-1}$
		holds on the open set $f(\Omega)$.
		Hence $\delta_{\Omega'}$ is differentiable at $f(z)$.
		By \eqref{eq:eikonal}, we have 
		$|\nabla\delta_\Omega(z)|
		=|\nabla\delta_{\Omega'}(f(z))|=1$.
		The chain rule and conformality therefore give
		\begin{equation}\label{eq:gradient-density}
			|f'(z)|=|\nabla(\delta_{\Omega'}\circ f)(z)|=e^{u(z)}
			|\nabla\delta_\Omega(z)+\delta_\Omega(z)\nabla u(z)|.
		\end{equation}
		Since $e^{u(z)}=|f'(z)|$, dividing
		\eqref{eq:gradient-density} by $e^{u(z)}$ and squaring gives
		\begin{equation}\label{eq:squared-gradient}
			1=1+2\delta_\Omega(z)\nabla\delta_\Omega(z)\cdot\nabla u(z)
			+\delta_\Omega(z)^2|\nabla u(z)|^2.
		\end{equation}
		Since $\delta_\Omega(z)>0$,
		\eqref{eq:squared-gradient} implies \eqref{eq:key-identity}. If $z\in G$, then \eqref{eq:u} gives
		$$ q(z)
		=-\frac{2}{P(z)}
		=-\frac{2\overline{P(z)}}{|P(z)|^2}
		=-\frac{2\nabla u(z)}{|\nabla u(z)|^2}.$$
		Thus \eqref{eq:key-identity} implies
		\eqref{eq:vector-identity}.
	\end{proof}
	
	\begin{remark}\label{rem:differential-identity-literature}
		After completing the first version of this manuscript, we learned that Graham and Witten had obtained the same differential equation as \eqref{eq:key-identity} in their study of geodesic defining functions for conformally compact Einstein metrics. For two such defining functions related by $\widehat r=e^\omega r$, their
		equation \cite{GW99}*{Lemma~2.2, equation~(2.8)} can be written as
		$$
		\langle dr,d\omega\rangle_{\bar g}=-\frac{r}{2}|d\omega|_{\bar g}^2,
		$$
		where $\bar g$ is the compactified metric associated with $r$. In our setting, $\bar g$ is the Euclidean metric, and the corresponding functions are $r=\delta_\Omega$,
		$\widehat r=\delta_{\Omega'}\circ f$, and $\omega=\log|f'|$. We also note that both \eqref{eq:key-identity} and \cite{GW99}*{equation~(2.8)} follow
		by differentiating the relation between the two functions and expanding the squared norm using their respective eikonal identities.
		
		We are not aware of an earlier use of \eqref{eq:key-identity} in the study of quasihyperbolic isometries.
	\end{remark}
	
	\section{Classification of the maps and the domains}
	\label{sec:classification}
	In this section we prove Theorem~\ref{thm:main}. We divide its proof into several lemmas.
	
	\begin{lemma}\label{lem:affine-vector}
		Let $\Omega,\Omega'\subsetneq\C$ be domains, and let
		$f:\Omega\to\Omega'$ be a surjective holomorphic
		quasihyperbolic isometry. Suppose that $P=f''/f'\not\equiv0$.
		Then $P$ has no zeros in $\Omega$, and there exist
		$\sigma\in\mathbb R$ and $p\in\C\setminus\Omega$
		such that
		\begin{equation}\label{eq:global-vector}
			P(z)=-\frac{2}{(1+i\sigma)(z-p)},\quad
			q(z):=-\frac{2}{P(z)}=(1+i\sigma)(z-p)
			\quad(z\in\Omega).
		\end{equation}
	\end{lemma}
	
	\begin{proof}
		Let $G=\{P\ne0\}$ and $q=-2/P$ on $G$. Fix $x\in G$ and
		choose $r$, $\xi$, and $n$ as in Lemma~\ref{lem:nearest-segment}.
		There exists $t_0\in(0,r)$ such that $\xi+tn\in G$ for
		$t_0<t<r$. Equations~\eqref{eq:nearest-segment} and
		\eqref{eq:vector-identity} imply
		\begin{equation}\label{eq:normal-vector}
			\operatorname{Re}\bigl(q(\xi+tn)\overline n\bigr)=t
			\quad(t_0<t<r).
		\end{equation}
		Differentiating \eqref{eq:normal-vector} with respect to $t$ gives
		\begin{equation}\label{eq:normal-derivative}
			\operatorname{Re}\bigl(q'(\xi+tn)n\overline n\bigr)
			=\operatorname{Re}q'(\xi+tn)=1.
		\end{equation}
		Letting $t\uparrow r$ in \eqref{eq:normal-derivative}, we obtain
		\begin{equation}\label{eq:real-derivative}
			\operatorname{Re}q'(x)=1\quad(x\in G).
		\end{equation}
		Choose a disk $V\subset G$. Since $q'$ is holomorphic,
		\eqref{eq:real-derivative} implies that, for some
		$\sigma\in\mathbb R$ and $b\in\C$,
		\begin{equation}\label{eq:local-affine}
			q'(z)=1+i\sigma,\quad q(z)=(1+i\sigma)z+b\quad(z\in V).
		\end{equation}
		By \eqref{eq:local-affine}, the holomorphic function
		$((1+i\sigma)z+b)P(z)+2$ vanishes on $V$. Therefore by the identity theorem we have 
		$$((1+i\sigma)z+b)P(z)=-2\quad(z\in\Omega).$$
		This gives \eqref{eq:global-vector}. In particular, $P$ has no zeros and $p=-b/(1+i\sigma)\notin\Omega$.
	\end{proof}
	
	\begin{lemma}\label{lem:flow}
		Let $\Omega\subsetneq\C$ be a domain, let $p\in\C$, and let $\sigma\in\mathbb R$. Suppose
		\begin{equation}\label{eq:flow-hypothesis}
			\nabla\delta_\Omega(z)\cdot(1+i\sigma)(z-p)=\delta_\Omega(z)
			\quad\text{for a.e. }z\in\Omega.
		\end{equation}
		Then $p\in\partial\Omega$ and, for
		$\Phi_t(z)=p+e^{(1+i\sigma)t}(z-p)$,
		\begin{equation}\label{eq:flow-invariance}
			\Phi_t(\Omega)=\Omega,\quad
			\delta_\Omega(\Phi_t(z))=e^t\delta_\Omega(z)
			\quad(z\in\Omega,\ t\in\mathbb R).
		\end{equation}
	\end{lemma}
	
	\begin{proof}
		Fix $x\in\Omega$ and choose an open neighbourhood
		$V\ni x$ with $\overline V\Subset\Omega$.
		Since $\Phi_t\to\operatorname{id}$ uniformly on
		$\overline V$ as $t\to0$, there exists $h>0$ such that
		$\Phi_t(V)\subset\Omega$ for $|t|<h$.
		
		Let $N\subset\Omega$ be a Borel null set outside which
		$\delta_\Omega$ is differentiable and
		\eqref{eq:flow-hypothesis} holds.
		For each fixed $t$, the set $\Phi_{-t}(N)$ has
		Lebesgue measure zero, since $\Phi_{-t}$ is a
		Euclidean similarity. Therefore, by Fubini's theorem it follows that for a.e. $z\in V$,
		$\Phi_t(z)\notin N$ for a.e. $t\in(-h,h)$.
		The chain rule and \eqref{eq:flow-hypothesis} give
		\begin{equation}\label{eq:flow-derivative}
			\begin{split}
				\frac{\mathrm d}{\mathrm dt}
				\bigl(e^{-t}\delta_\Omega(\Phi_t(z))\bigr)
				&=e^{-t}\bigl(
				\nabla\delta_\Omega(\Phi_t(z))
				\cdot(1+i\sigma)(\Phi_t(z)-p)
				-\delta_\Omega(\Phi_t(z))\bigr)\\
				&=0\quad\text{for a.e. }t\in(-h,h).
			\end{split}
		\end{equation}
		Since $\delta_\Omega$ is
		$1$-Lipschitz and $\Phi_t$ is smooth, for each fixed $z\in V$, the function
		$t\mapsto e^{-t}\delta_\Omega(\Phi_t(z))$
		is locally Lipschitz. Thus \eqref{eq:flow-derivative} implies that this
		function is constant for a.e. $z\in V$.
		Continuity in $z$ then gives
		\begin{equation}\label{eq:local-flow-scaling}
			\delta_\Omega(\Phi_t(z))=e^t\delta_\Omega(z)
			\quad(z\in V,\ |t|<h).
		\end{equation}
		Fix $z\in\Omega$ and let $I_z$ be the component of
		$\{t\in\mathbb R:\Phi_t(z)\in\Omega\}$ containing $0$.
		For every $\tau\in I_z$, apply
		\eqref{eq:local-flow-scaling} at $\Phi_\tau(z)$.
		Using $\Phi_s\circ\Phi_\tau=\Phi_{\tau+s}$, we obtain
		\[
		\delta_\Omega(\Phi_{\tau+s}(z))
		=e^s\delta_\Omega(\Phi_\tau(z))
		\quad\text{for all sufficiently small }|s|.
		\]
		Hence $t\mapsto e^{-t}\delta_\Omega(\Phi_t(z))$
		is locally constant on the interval $I_z$.
		It is therefore constant on $I_z$, and evaluation
		at $t=0$ gives
		\begin{equation}\label{eq:maximal-flow-scaling}
			\delta_\Omega(\Phi_t(z))=e^t\delta_\Omega(z)
			\quad(t\in I_z).
		\end{equation}
		
		If $I_z$ had a finite endpoint $T$, continuity of
		the distance to the closed complement and
		\eqref{eq:maximal-flow-scaling} would imply
		\begin{equation}\label{eq:no-finite-exit}
			\dist(\Phi_T(z),\C\setminus\Omega)
			=\lim_{t\to T,\ t\in I_z}\delta_\Omega(\Phi_t(z))
			=e^T\delta_\Omega(z)>0.
		\end{equation}
		Thus $\Phi_T(z)\in\Omega$. Since $\Omega$ is open,
		this contradicts the definition of $I_z$.
		Hence $I_z=\mathbb R$.
		Since $\Phi_{-t}=\Phi_t^{-1}$,
		\eqref{eq:maximal-flow-scaling} proves
		\eqref{eq:flow-invariance}.
		
		Finally, fix $z\in\Omega$. Since
		$\Phi_t(z)\to p$ as $t\to-\infty$,
		we have $p\in\overline\Omega$. Moreover,
		\[
		\dist(p,\C\setminus\Omega)
		=\lim_{t\to-\infty}\delta_\Omega(\Phi_t(z))
		=\lim_{t\to-\infty}e^t\delta_\Omega(z)=0.
		\]
		Therefore $p\in\partial\Omega$.
	\end{proof}
	
	\begin{lemma}\label{lem:spiral-domains}
		Let $p\in\C$, $\sigma\in\mathbb{R}$, and
		$\Omega\subsetneq\C$ be a domain. If $\Omega$ satisfies $\Phi_t(\Omega)=\Omega$ for every $t\in \mathbb{R}$,
		then either $\Omega=\C\setminus\{p\}$ or
		\begin{equation}\label{eq:domain-classification}
			\Omega=p+D_{\sigma,\theta,\alpha}
			\quad\text{for some }\theta\in\mathbb R,\quad0<\alpha\leq2\pi.
		\end{equation}
	\end{lemma}
	
	\begin{proof}
		We first prove that $p$ cannot be in $\Omega$.
		If $p\in\Omega$, choose $r>0$ with $B(p,r)\subset\Omega$.
		Since $\Phi_t(\Omega)=\Omega$ for every $t\in\mathbb{R}$,
		it follows that $B(p,e^tr)\subset\Omega$ for every
		$t\in\mathbb{R}$, contrary to $\Omega\ne\C$.
		Thus $p\notin\Omega$. After translation, assume $p=0$.
		Define $S^1=\{\zeta\in\C:|\zeta|=1\}$ and the homeomorphism $H:\C\setminus\{0\}\to \mathbb{R}\times S^1$ by
		\begin{equation}\label{eq:untwisting}
			\begin{split}
				H(z)&=\left(\log|z|,
				\frac{z}{|z|}e^{-i\sigma\log|z|}\right),\\
				H^{-1}(s,\zeta)&=e^{(1+i\sigma)s}\zeta.
			\end{split}
		\end{equation}
		The invariance $\Phi_t(\Omega)=\Omega$ for every
		$t\in\mathbb{R}$, together with \eqref{eq:untwisting}, gives
		\begin{equation}\label{eq:product-domain}
			H(\Omega)=\mathbb{R}\times V,
		\end{equation}
		where $V$ is a nonempty open subset of $S^1$.
		Since $\Omega$ is connected, so is $H(\Omega)$,
		and hence its projection $V$ onto $S^1$. If $V=S^1$, then $\Omega=\C\setminus\{0\}$.
		Otherwise,
		\begin{equation}\label{eq:angular-interval}
			V=\{e^{i\beta}:\theta<\beta<\theta+\alpha\}
			\quad\text{for some }\theta\in\mathbb{R}
			\text{ and }0<\alpha\leq2\pi.
		\end{equation}
		By equations~\eqref{eq:untwisting}, \eqref{eq:product-domain}, and \eqref{eq:angular-interval} it follows that $\Omega=D_{\sigma,\theta,\alpha}$.
	\end{proof}
	
	\subsection{Proof of the main result}\label{subsec:proof-of-main-result}
	
	\begin{proof}[Proof of Theorem~\ref{thm:main}]
		By \cite{MO86}*{Theorem~2.6}, after postcomposition
		with a reflection if necessary, we may assume that
		$f$ is holomorphic. We consider two cases: $P\equiv0$ and $P\not\equiv0$. If $P\equiv0$, then $f(z)=az+b$ with $a\ne0$. Thus, we have Theorem~\ref{thm:main-similarity}.
		
		If $P\not\equiv0$, we first provide the domain classification.
		Since $P\not\equiv0$, $G\neq\varnothing$, and by
		Lemma~\ref{lem:affine-vector}, there exist
		$\sigma\in\mathbb{R}$ and $p\notin\Omega$ such that
		\eqref{eq:global-vector} holds.
		In particular, \eqref{eq:global-vector} gives $G=\Omega$.
		Since $\delta_\Omega$ is $1$-Lipschitz, it is differentiable
		a.e. in $\Omega$. Thus, by Lemma~\ref{lem:identity},
		\eqref{eq:flow-hypothesis} holds for a.e. $z\in\Omega$.
		Hence, by Lemma~\ref{lem:flow}, we have
		$p\in\partial\Omega$ and \eqref{eq:flow-invariance}.
		Finally, Lemma~\ref{lem:spiral-domains} gives the
		classification of $\Omega$, i.e., $\Omega$ is either
		$\mathbb{C}\setminus\{p\}$ or $p+D_{\sigma,\theta,\alpha}$
		for some $\theta\in\mathbb{R}$ and $0<\alpha\leq2\pi$.
		
		Now we study two cases: $\sigma=0$ and $\sigma\neq 0$. If $\sigma=0$, \eqref{eq:global-vector} gives
		\begin{equation}\label{eq:integrate-reciprocal}
			\bigl(f'(z)(z-p)^2\bigr)'=0,
			\quad f'(z)=\frac{C}{(z-p)^2},\quad
			f(z)=B-\frac{C}{z-p}
		\end{equation}
		for some $C\ne0$. This proves Theorem~\ref{thm:main-reciprocal}.
		
		Assume $\sigma\ne0$. We show that $\Omega\neq \mathbb{C}\setminus \{p\}$. If $\Omega=\C\setminus\{p\}$,
		then, for the positively oriented circle
		$\gamma_r=\{z:|z-p|=r\}$,
		\begin{equation}\label{eq:period-obstruction}
			\mathbb Z\ni
			\frac1{2\pi i}\int_{\gamma_r}
			\frac{f''(z)}{f'(z)}\,\mathrm dz
			=-\frac2{1+i\sigma}\notin\mathbb R,
		\end{equation}
		a contradiction. Here the integral divided by
		$2\pi i$ is the winding number of
		$f'\circ\gamma_r$ around $0$.
		Hence $\Omega=p+D_{\sigma,\theta,\alpha}$. For simplicity, write $L=L_{\sigma,\theta,\alpha}$ and take
		$F(z)=\exp(\lambda_\sigma L(z-p))$, then
		\begin{equation}\label{eq:integrate-power}
			\frac{F''(z)}{F'(z)}
			=\frac{\lambda_\sigma-1}{z-p}
			=-\frac2{(1+i\sigma)(z-p)}
			=\frac{f''(z)}{f'(z)}.
		\end{equation}
		Thus $(f'/F')'=0$ and $f=AF+B$ for some $A\ne0$.
		This proves \eqref{eq:main-power}.
		
		Conversely since Euclidean similarities are quasihyperbolic isometries, it suffices to prove the converse in the
		sector cases with $p=B=0$ and $A=1$. Write
		\begin{equation}\label{eq:normalized-power}
			D=D_{\sigma,\theta,\alpha},\quad
			U=U_{\sigma,\theta,\alpha},\quad L=(\exp|_U)^{-1},\quad
			F(z)=\exp(\lambda_\sigma L(z)),
		\end{equation}
		where $\sigma=0$ is allowed. For $z\in D$, write $w=L(z)=s+it$. It follows that
		\begin{equation}\label{eq:strip-reflection}
			\operatorname{Im}(\lambda_\sigma w)
			+\sigma\operatorname{Re}(\lambda_\sigma w)
			=-(t-\sigma s),\quad \lambda_\sigma U=U^*.
		\end{equation}
		The exponential is injective on $U$ and $U^*$, so
		$F:D\to D^*$ is a conformal bijection. Define
		\begin{equation}\label{eq:power-similarity}
			\tau=-\frac{2(s+\sigma t)}{1+\sigma^2},\quad
			\lambda_\sigma w-\overline w=(1-i\sigma)\tau,
			\quad F(z)=e^{(1-i\sigma)\tau}\overline z.
		\end{equation}
		Since $U+(1+i\sigma)v=U$ for $v\in\mathbb R$,
		\begin{equation}\label{eq:reflected-scaling}
			e^{(1-i\sigma)v}D^*=D^*,\quad
			\delta_{D^*}(e^{(1-i\sigma)v}\zeta)
			=e^v\delta_{D^*}(\zeta)
			\quad(\zeta\in D^*,\ v\in\mathbb R).
		\end{equation}
		Equations~\eqref{eq:power-similarity} and
		\eqref{eq:reflected-scaling} imply
		\begin{equation}\label{eq:converse-density}
			\begin{split}
				\delta_{D^*}(F(z))&=e^\tau\delta_{D^*}(\overline z)
				=e^\tau\delta_D(z),\\
				|F'(z)|&=\left|\lambda_\sigma\frac{F(z)}z\right|
				=e^{\operatorname{Re}(\lambda_\sigma w)-s}=e^\tau.
			\end{split}
		\end{equation}
		Lemma~\ref{lem:density} applies to \eqref{eq:converse-density}.
		For the punctured plane, we have
		$$\delta_{\C\setminus\{0\}}(1/z)=\frac1{|z|}
		=\left|\left(\frac1z\right)'\right|
		\delta_{\C\setminus\{0\}}(z),$$
		so Lemma~\ref{lem:density} proves the remaining case.
		Finally, substituting \eqref{eq:integrate-power} into
		\eqref{eq:schwarzian-definition} gives
		\eqref{eq:main-schwarzian}.
	\end{proof}
	
	\subsubsection{Quasihyperbolic isometry classification
		in dimensions \texorpdfstring{$n\geq3$}{n >= 3}}
	\label{sec:higher}
	\begin{proof}[Proof of Theorem~\ref{Thm:higher-dimension-classification}]
		By \cite{MO86}*{Theorem~2.6} and Liouville's theorem,
		$f$ is the restriction of a M\"obius transformation.
		If $f$ is not a similarity, its standard decomposition gives
		\[
		f(x)=b+cQJ(x-p),\quad
		J(y)=\frac{y}{|y|^2},
		\]
		where $p\notin\Omega$, $b\in\rn$, $c>0$, and $Q\in O(n)$.
		Since similarities are quasihyperbolic isometries,
		$J:\Omega-p\to J(\Omega-p)$ is a quasihyperbolic isometry.
		Thus \cite{R19}*{Theorem~7.1.5} gives
		$p\in\partial\Omega$ and
		$p+t(\Omega-p)=\Omega$ for every $t>0$.
		
		Conversely, suppose $f$ has the form
		\eqref{eq:higher-classification} and write $D=\Omega-p$.
		The identities $tD=D$ and $\delta_D(tx)=t\delta_D(x)$
		for $t>0$ imply that $J(D)=D$ and
		\begin{equation}\label{eq:higher-inversion-density}
			\delta_D(J(x))=\frac{\delta_D(x)}{|x|^2},
			\quad
			|DJ(x)v|=\frac{|v|}{|x|^2}
			\quad(x\in D,\ v\in\rn).
		\end{equation}
		By \eqref{eq:higher-inversion-density}, $J$ preserves
		quasihyperbolic lengths. Since $J^{-1}=J$, it is a
		quasihyperbolic isometry of $D$.
		Composition with similarities proves the converse,
		and $J(D)=D$ gives \eqref{eq:higher-image}.
		The similarity case follows directly from the definition.
	\end{proof}
	
	\begin{proof}[Proof of Corollary~\ref{cor:inradius-all-dimensions}]
		For $n=2$, it follows from Corollary~\ref{cor:inradius}.
		For $n\geq3$, replacing $f$ by $f^{-1}$ if necessary, assume
		$R_\Omega<\infty$. If $f$ is not a similarity,
		\eqref{eq:higher-classification} gives
		$$\delta_\Omega(p+t(x-p))=t\delta_\Omega(x)
		\longrightarrow\infty\quad(t\to\infty)$$
		for every $x\in\Omega$, a contradiction.
	\end{proof}

	\section{Consequences and examples}
	\label{sec:consequences}
	\begin{proof}[Proof of Corollary~\ref{cor:inradius}]
		Suppose $f$ is not a similarity. Theorem~\ref{thm:main} and
		\eqref{eq:flow-invariance} give, for any fixed $z\in\Omega$,
		\begin{equation}\label{eq:infinite-inradius}
			R_\Omega\geq\delta_\Omega(\Phi_t(z))
			=e^t\delta_\Omega(z)\longrightarrow\infty
			\quad(t\to\infty).
		\end{equation}
		The inverse $f^{-1}$ is also a nonsimilarity isometry, so the same
		argument gives $R_{\Omega'}=\infty$.
	\end{proof}
	
	\begin{proof}[Proof of Corollary~\ref{cor:connectivity}]
		Each $D_{\sigma,\theta,\alpha}$ is conformally
		equivalent to the strip $U_{\sigma,\theta,\alpha}$
		and hence is simply connected. The conclusion
		follows from Theorem~\ref{thm:main}.
	\end{proof}
	
	\begin{proof}[Proof of Corollary~\ref{cor:c1}]
		Suppose that $f$ is not a similarity.
		By Theorem~\ref{thm:main}, $\Omega$ is a punctured
		plane or $p+D_{\sigma,\theta,\alpha}$.
		The punctured plane has an isolated boundary point,
		while for $\sigma\ne0$ the boundary has no tangent
		line at $p$. For $\sigma=0$, the sector boundary
		is locally a $C^1$ graph at $p$ if and only if
		$\alpha=\pi$. Thus $\Omega$ is a half-plane,
		a contradiction.
	\end{proof}
	
	\begin{proof}[Proof of Corollary~\ref{cor:similar-domains}]
		In the sector cases, \eqref{eq:strip-reflection} and
		Theorem~\ref{thm:main} give, for the holomorphic normalization,
		\begin{equation}\label{eq:similar-domain-image}
			\Omega'=B+A(\Omega-p)^*.
		\end{equation}
		The map $z\mapsto B+A\overline{z-p}$ is a Euclidean similarity
		realizing \eqref{eq:similar-domain-image}.
		The punctured-plane case follows from
		$f(\C\setminus\{p\})=\C\setminus\{B\}$.
		The similarity case and the target reflection preserve the conclusion.
	\end{proof}
	
	\begin{proof}[Proof of Corollary~\ref{cor:self-rigidity}]
		It is clear that each rigidity property implies the corresponding self-rigidity property. Conversely, let
		$f:(\Omega,k_\Omega)\to(f(\Omega),k_{f(\Omega)})$
		be a quasihyperbolic isometry, where $f(\Omega)$
		is a proper planar domain. By Corollary~\ref{cor:similar-domains}, choose a
		Euclidean similarity $S$ with $S(\Omega)=f(\Omega)$.
		Then $g=S^{-1}\circ f$ is a surjective
		quasihyperbolic self-isometry of $\Omega$.
		If $\Omega$ is self-similar rigid, then
		$f=S\circ g$ is a similarity. If $\Omega$ is
		self-M\"obius rigid, then $f=S\circ g$ is M\"obius.
	\end{proof}
	
	\begin{example}\label{ex:spiral}
		Let
		\begin{equation}\label{eq:spiral-example}
			\begin{split}
				U&=\{s+it:0<t-s<\pi\},\quad \Omega=\exp U,\\
				L&=(\exp|_U)^{-1},\quad f(z)=\exp(iL(z)).
			\end{split}
		\end{equation}
		Since $iU=U^*$, the map $f:\Omega\to\Omega^*$ is a conformal
		bijection. For $L(z)=s+it$ and $\tau=-s-t$,
		\begin{equation}\label{eq:explicit-spiral-density}
			\begin{split}
				f(z)&=e^{(1-i)\tau}\overline z,\quad
				f'(z)=\frac{if(z)}z,\quad |f'(z)|=e^\tau,\\
				\delta_{\Omega^*}(f(z))
				&=e^\tau\delta_{\Omega^*}(\overline z)
				=|f'(z)|\delta_\Omega(z).
			\end{split}
		\end{equation}
		Lemma~\ref{lem:density} and
		\eqref{eq:explicit-spiral-density} show that $f$ is a
		quasihyperbolic isometry. However,
		\begin{equation}\label{eq:explicit-schwarzian}
			\frac{f''}{f'}=\frac{i-1}{z},\quad
			S_f(z)=-\frac{i-1}{z^2}-\frac{(i-1)^2}{2z^2}
			=\frac1{z^2}\ne0.
		\end{equation}
		Thus $f$ is not M\"obius. The same construction also
		gives an explicit failure of self-M\"obius rigidity.
		Indeed, $ g(z)=\overline{f(z)}$ is an anticonformal quasihyperbolic self-isometry
		of $\Omega$. Moreover,
		\begin{equation}\label{eq:spiral-involution}
			L(g(z))=-i\overline{L(z)},\quad
			g\circ g=\operatorname{id}_\Omega.
		\end{equation}
		If $g$ were M\"obius, then so would be
		$f=\overline g$, contradicting
		\eqref{eq:explicit-schwarzian}. Hence, $g$ is not M\"obius. By \eqref{eq:infinite-inradius},
		$R_\Omega=R_{\Omega^*}=\infty$.
	\end{example}
	
	\section{Proof of Theorem~\ref{thm:planar-pgh}}
	\label{sec:pgh-proof}
	We use pointed Gromov--Hausdorff convergence in the
	sense of \cite{HRS20}*{Section~2.3}. We first strengthen the similarity conclusion for
	quasihyperbolically isometric domains by preserving a prescribed base point.
	
	\begin{lemma}\label{lem:pointed-similarity}
		Let $\Omega,\Omega'\subsetneq\C$ be domains, let
		$f:(\Omega,k_\Omega)\to(\Omega',k_{\Omega'})$ be a
		surjective isometry, and let $a\in\Omega$. Then there exist a Euclidean similarity $S:\C\to\C$
		and $c>0$ such that
		$|S(z)-S(w)|=c|z-w|$ for all $z,w\in\C$, and
		\begin{equation}\label{eq:pointed-similarity}
			S(\Omega)=\Omega',\quad S(a)=f(a),\quad
			c=\frac{\delta_{\Omega'}(f(a))}{\delta_\Omega(a)}.
		\end{equation}
	\end{lemma}
	
	\begin{proof}
		The assertion is preserved under pre- and postcomposition
		by Euclidean similarities. If $f$ is a similarity, take
		$S=f$. Otherwise, by Theorem~\ref{thm:main}, we may assume
		that $f$ is holomorphic and that $p=B=0$ and $A=1$. If $\Omega=\C\setminus\{0\}$ and $f(z)=1/z$, take
		$S(z)=\overline z/|a|^2$. In the remaining cases,
		$\Omega=D_{\sigma,\theta,\alpha}$ and $f$ has the form
		\eqref{eq:normalized-power}, with $\sigma=0$ allowed.
		Write
		$$w=L_{\sigma,\theta,\alpha}(a),\quad
		t_0=-2\operatorname{Re}\frac{w}{1+i\sigma},\quad 
		S(z)=e^{(1-i\sigma)t_0}\overline z
		\quad(z\in\C).$$
		Equations~\eqref{eq:power-similarity} and
		\eqref{eq:reflected-scaling} give
		$S(a)=f(a)$ and $S(\Omega)=f(\Omega)$.
		Finally, if $c>0$ is the scale factor of $S$, then
		$\delta_{\Omega'}(f(a))=c\delta_\Omega(a)$.
	\end{proof}
	
	In particular, for proper planar domains $D,U$
	containing $0$, Lemma~\ref{lem:pointed-similarity} gives
	\begin{equation}\label{eq:pointed-orthogonal-equivalence}
		\left.
		\begin{aligned}
			(D,k_D;0)&\cong(U,k_U;0),\\
			\delta_D(0)&=\delta_U(0)
		\end{aligned}
		\right\}
		\quad\Longrightarrow\quad
		Q(D)=U
		\quad\text{for some }Q\in O(2),
	\end{equation}
	where $\cong$ denotes pointed isometry.
	
	\begin{proof}[Proof of Theorem~\ref{thm:planar-pgh}]
		Write $r=\delta_\Omega(a)$, $r_j=\delta_{\Omega_j}(a_j)$,
		\begin{equation}\label{eq:pgh-normalization}
			T(z)=\frac{z-a}{r}, \quad T_j(z)=\frac{z-a_j}{r_j}, \quad U=T(\Omega), \quad U_j=T_j(\Omega_j).
		\end{equation}
		Then
		$\delta_U(0)=\delta_{U_j}(0)=1$. Suppose \eqref{eq:planar-pgh-convergence} holds.
		Similarity invariance gives
		\[
		(U_j,k_{U_j};0)
		\xrightarrow{\mathrm{pGH}}
		(U,k_U;0).
		\]
		In the proof of \cite{HRS20}*{Theorem~5.8},
		the assumption $n\geq3$ is used only to deduce that
		two pointed quasihyperbolically isometric domains
		with equal boundary distances at the origin are
		related by an orthogonal transformation.
		Equation~\eqref{eq:pointed-orthogonal-equivalence}
		provides this conclusion in dimension two. Therefore, following the argument in the proof of
		\cite{HRS20}*{Theorem~5.8} for $U_j$ and $U$ with
		\eqref{eq:pointed-orthogonal-equivalence},
		we obtain $Q_j\in O(2)$ such that
		\begin{equation}\label{eq:pgh-normalized-convergence}
			d_{Q_j(U_j)}\longrightarrow\delta_U
			\quad\text{uniformly on bounded subsets of }U.
		\end{equation}
		For the similarities in
		\eqref{eq:planar-pgh-similarities}, we have
		\begin{equation}\label{eq:pgh-scaling-identities}
			\begin{aligned}
				S_j&=T^{-1}\circ Q_j\circ T_j,\\
				d_{S_j(\Omega_j)}(z)
				&=r\,d_{Q_j(U_j)}(T(z)),\\
				\delta_\Omega(z)
				&=r\,\delta_U(T(z))
				\quad(z\in\Omega).
			\end{aligned}
		\end{equation}
		Equations~\eqref{eq:pgh-normalized-convergence}
		and \eqref{eq:pgh-scaling-identities} imply
		\eqref{eq:planar-pgh-uniform}.
		
		Conversely, since every closed quasihyperbolic ball $K$ in $\Omega$ is compactly contained in $\Omega$, \eqref{eq:planar-pgh-uniform} gives
		$$ K\subset V_j
		\quad\text{for all sufficiently large }j\, \text{ and }\,
		\frac{\delta_\Omega}{d_{V_j}}
		\longrightarrow1
		\quad\text{uniformly on }K,$$
		where $V_j=S_j(\Omega_j)$. By \cite{HRS20}*{Theorem~4.4}, $(V_j,k_{V_j};a)
		\xrightarrow{\mathrm{pGH}}
		(\Omega,k_\Omega;a)$. Each $S_j$ is a quasihyperbolic isometry with
		$S_j(a_j)=a$, so
		\eqref{eq:planar-pgh-convergence} follows. Finally, $d_{V_j}$ and $\delta_\Omega$ are $1$-Lipschitz on $\Omega$. A finite-net argument on each bounded subset of $\Omega$ shows that pointwise convergence is equivalent to \eqref{eq:planar-pgh-uniform}; compare
		\cite{HRS20}*{the proof of Proposition~3.8 and Remark~3.9(a)}.
	\end{proof}
	
	\section*{Declaration of generative AI use}
	The author used ChatGPT (OpenAI, GPT-6 Astra) in preparing
	the proof of Theorem~\ref{thm:planar-pgh}. Specifically, interactions with ChatGPT helped the author identify
	\eqref{eq:pointed-orthogonal-equivalence} as the planar
	replacement for the step in the argument of Herron, Richard,
	and Snipes that requires $n\geq3$. Comparing their proof
	with the planar argument with the help of ChatGPT allowed
	the author to shorten the proof, as explained in the final
	paragraph of the proof of Theorem~\ref{thm:planar-pgh}. A literature search with ChatGPT for earlier occurrences of
	\eqref{eq:key-identity} helped the author locate the equation of Graham and Witten discussed in
	Remark~\ref{rem:differential-identity-literature}. The author also used ChatGPT to assist with language editing and improving the exposition of the manuscript. Apart from the assistance described above, all mathematical statements and all other proofs were developed independently by the author. The author takes full responsibility for the content of the manuscript.


\begin{bibdiv}
		\begin{biblist}
			
			\bib{BB03}{article}{
				author={Balogh, Zolt\'an M.},
				author={Buckley, Stephen M.},
				title={Geometric characterizations of Gromov hyperbolicity},
				journal={Invent. Math.},
				volume={153},
				number={2},
				date={2003},
				pages={261--301},
			}
			
			\bib{BHK01}{book}{
				author={Bonk, Mario},
				author={Heinonen, Juha},
				author={Koskela, Pekka},
				title={Uniformizing Gromov hyperbolic spaces},
				series={Ast\'erisque},
				volume={270},
				publisher={Soci\'et\'e Math\'ematique de France},
				date={2001},
			}
			
			\bib{BH20}{article}{
				author={Buckley, Stephen M.},
				author={Herron, David A.},
				title={Quasihyperbolic geodesics are hyperbolic quasi-geodesics},
				journal={J. Eur. Math. Soc. (JEMS)},
				volume={22},
				number={6},
				date={2020},
				pages={1917--1970},
			}
			
			\bib{CF79}{article}{
				author={Caffarelli, Luis A.},
				author={Friedman, Avner},
				title={The free boundary for elastic-plastic torsion problems},
				journal={Trans. Amer. Math. Soc.},
				volume={252},
				date={1979},
				pages={65--97},
			}
			
			\bib{CCM97}{article}{
				author={Choi, Hyeong In},
				author={Choi, Sung Woo},
				author={Moon, Hwan Pyo},
				title={Mathematical theory of medial axis transform},
				journal={Pacific J. Math.},
				volume={181},
				number={1},
				date={1997},
				pages={57--88},
			}
			
			\bib{Dam03}{article}{
				author={Damon, James},
				title={Smoothness and geometry of boundaries associated to
					skeletal structures I: sufficient conditions for smoothness},
				journal={Ann. Inst. Fourier (Grenoble)},
				volume={53},
				number={6},
				date={2003},
				pages={1941--1985},
			}
			
			\bib{GO79}{article}{
				author={Gehring, F. W.},
				author={Osgood, Brad G.},
				title={Uniform domains and the quasi-hyperbolic metric},
				journal={J. Anal. Math.},
				volume={36},
				date={1979},
				pages={50--74},
			}
			
			\bib{GP76}{article}{
				author={Gehring, F. W.},
				author={Palka, B. P.},
				title={Quasiconformally homogeneous domains},
				journal={J. Anal. Math.},
				volume={30},
				date={1976},
				pages={172--199},
			}
			
			\bib{GW99}{article}{
				author={Graham, C. Robin},
				author={Witten, Edward},
				title={Conformal anomaly of submanifold observables
					in {AdS/CFT} correspondence},
				journal={Nuclear Phys. B},
				volume={546},
				number={1--2},
				date={1999},
				pages={52--64},
			}
			
			\bib{H07}{article}{
				author={H\"ast\"o, Peter},
				title={Isometries of the quasihyperbolic metric},
				journal={Pacific J. Math.},
				volume={230},
				number={2},
				date={2007},
				pages={315--326},
			}
			
			\bib{HI05}{article}{
				author={H\"ast\"o, Peter},
				author={Ibragimov, Zair},
				title={Apollonian isometries of planar domains are
					M\"obius mappings},
				journal={J. Geom. Anal.},
				volume={15},
				number={2},
				date={2005},
				pages={229--237},
			}
			
			\bib{HIL06}{article}{
				author={H\"ast\"o, Peter},
				author={Ibragimov, Zair},
				author={Lind\'en, Henri},
				title={Isometries of relative metrics},
				journal={Comput. Methods Funct. Theory},
				volume={6},
				number={1},
				date={2006},
				pages={15--28},
			}
			
			\bib{HIMPS07}{article}{
				author={H\"ast\"o, Peter},
				author={Ibragimov, Zair},
				author={Minda, David},
				author={Ponnusamy, Saminathan},
				author={Sahoo, Swadesh},
				title={Isometries of some hyperbolic-type path metrics,
					and the hyperbolic medial axis},
				book={
					title={In the tradition of Ahlfors--Bers. IV},
					series={Contemp. Math.},
					volume={432},
					publisher={Amer. Math. Soc.},
					address={Providence, RI},
				},
				date={2007},
				pages={63--74},
			}
			
			\bib{Hei62}{article}{
				author={Heins, Maurice},
				title={On a class of conformal metrics},
				journal={Nagoya Math. J.},
				volume={21},
				date={1962},
				pages={1--60},
			}
			
			\bib{Her21}{article}{
				author={Herron, David A.},
				title={Quasihyperbolic metric universal covers},
				journal={Rev. Mat. Iberoam.},
				volume={37},
				number={4},
				date={2021},
				pages={1277--1306},
			}
			
			\bib{HIM08}{article}{
				author={Herron, David A.},
				author={Ibragimov, Zair},
				author={Minda, David},
				title={Geodesics and curvature of M\"obius invariant metrics},
				journal={Rocky Mountain J. Math.},
				volume={38},
				number={3},
				date={2008},
				pages={891--921},
			}
			
			\bib{HRS20}{article}{
				author={Herron, David A.},
				author={Richard, Abigail},
				author={Snipes, Marie A.},
				title={Chordal Hausdorff convergence and quasihyperbolic distance},
				journal={Anal. Geom. Metr. Spaces},
				volume={8},
				date={2020},
				pages={36--67},
			}
			
			\bib{IP26}{misc}{
				author={Ikonen, Toni},
				author={Pandey, Abhishek},
				title={Quasihyperbolic domains are $\mathrm{CAT}(2)$},
				date={2026},
				note={Preprint, arXiv:2609.24610},
			}
			
			\bib{Jon80}{article}{
				author={Jones, Peter W.},
				title={Extension theorems for {BMO}},
				journal={Indiana Univ. Math. J.},
				volume={29},
				number={1},
				date={1980},
				pages={41--66},
			}
			
			\bib{JS00}{article}{
				author={Jones, Peter W.},
				author={Smirnov, Stanislav K.},
				title={Removability theorems for Sobolev functions
					and quasiconformal maps},
				journal={Ark. Mat.},
				volume={38},
				number={2},
				date={2000},
				pages={263--279},
			}
			
			\bib{Kos98}{article}{
				author={Koskela, Pekka},
				title={Old and new on the quasihyperbolic metric},
				book={
					title={Quasiconformal mappings and analysis},
					editor={Duren, Peter},
					editor={Heinonen, Juha},
					editor={Osgood, Brad},
					editor={Palka, Bruce},
					publisher={Springer},
					address={New York},
				},
				date={1998},
				pages={205--219},
			}
			
			\bib{MM03}{article}{
				author={Mantegazza, Carlo},
				author={Mennucci, Andrea Carlo},
				title={Hamilton--Jacobi equations and distance functions
					on Riemannian manifolds},
				journal={Appl. Math. Optim.},
				volume={47},
				number={1},
				date={2003},
				pages={1--25},
			}
			
			\bib{Mar85}{article}{
				author={Martin, Gaven J.},
				title={Quasiconformal and bi-Lipschitz homeomorphisms,
					uniform domains and the quasihyperbolic metric},
				journal={Trans. Amer. Math. Soc.},
				volume={292},
				number={1},
				date={1985},
				pages={169--191},
			}
			
			\bib{MO86}{article}{
				author={Martin, Gaven J.},
				author={Osgood, Brad G.},
				title={The quasihyperbolic metric and associated estimates
					on the hyperbolic metric},
				journal={J. Anal. Math.},
				volume={47},
				date={1986},
				pages={37--53},
			}
			
			
			
			
			\bib{NY20}{article}{
				author={Ntalampekos, Dimitrios},
				author={Younsi, Malik},
				title={Rigidity theorems for circle domains},
				journal={Invent. Math.},
				volume={220},
				number={1},
				date={2020},
				pages={129--183},
			}
			
			\bib{Osg98}{article}{
				author={Osgood, Brad},
				title={Old and new on the Schwarzian derivative},
				book={
					title={Quasiconformal mappings and analysis},
					editor={Duren, Peter},
					editor={Heinonen, Juha},
					editor={Osgood, Brad},
					editor={Palka, Bruce},
					publisher={Springer},
					address={New York},
				},
				date={1998},
				pages={275--308},
			}
			
			\bib{R19}{thesis}{
				author={Richard, Abigail},
				title={Quasihyperbolic distance, pointed Gromov--Hausdorff
					distance, and bounded uniform convergence},
				type={Ph.D. thesis},
				organization={University of Cincinnati},
				date={2019},
			}
			
			\bib{Vai09}{article}{
				author={V\"ais\"al\"a, Jussi},
				title={Quasihyperbolic geometry of planar domains},
				journal={Ann. Acad. Sci. Fenn. Math.},
				volume={34},
				number={2},
				date={2009},
				pages={447--473},
			}
			
		\end{biblist}
	\end{bibdiv}
\end{document}